\documentclass[review,authoryear]{elsarticle}
\usepackage{amssymb,amsfonts,amsthm,amsmath,graphicx,mathrsfs,url,bm,bbm}

\newcommand{\beq}{\begin{equation}}
\newcommand{\eeq}{\end{equation}}

\newtheorem{theorem}{Theorem}

\newtheorem{corollary}[theorem]{Corollary}

\def\bvec#1{\mbox{\boldmath $#1$}}

\def\bone{{\bvec 1}}

\def\V{{\rm Var}}

\usepackage{mathrsfs,url,color}
\def\convd{{\stackrel{d}{\rightarrow}}}

\def\equald{{\stackrel{d}{=}}}
\def\ignore#1{}
\def\Ref#1{(\ref{#1})}

\newcommand{\ScB}{\mathscr{B}}
\newcommand{\ScF}{\mathscr{F}}
\newcommand{\ScS}{\mathscr{S}}
\newcommand{\ScM}{\mathscr{M}}
\newcommand{\ScN}{\mathscr{N}}
\newcommand{\ScX}{\mathscr{X}}
\newcommand{\prob}{\mathbb{P}}
\newcommand{\mean}{\mathbb{E}}

\newcommand{\Z}{\mathbb{Z}}
\newcommand{\N}{\mathbb{N}}

\def\bonee{{\bf\scriptsize 1}}

\def\ime{{\mbox{\sl\scriptsize i}}}
\def\imee{{\mbox{\sl\tiny i}}}
\def\Pn{{\rm Pn}}

\begin{document}

\begin{frontmatter}

	\title{A {marked renewal process model} for the size of a honey bee {colony}}

	\author[um]{Aihua Xia\corref{cor1}\fnref{fn1}}
	\ead{aihuaxia@unimelb.edu.au}
	\author[um]{Richard M. Huggins}
	\author[mb]{Martine J. Barons\fnref{fn2}}
	\author[um]{Louis Guillot}
	\address[um]{School of Mathematics and Statistics, The University of Melbourne, Parkville, Victoria, 3010
Australia}
	\address[mb]{Department of Statistics, The University of Warwick, Coventry
CV4 7AL
United Kingdom}

	\cortext[cor1]{Corresponding author}

	\fntext[fn1]{Research supported by ARC Discovery DP150101459.}
	\fntext[fn2]{Research supported by EPSRC grant EP/K039628/1.}

		\begin{abstract}
	Many areas of agriculture rely on honey bees to provide pollination services and any decline in honey bee numbers can impact on global food security. 
	In order to understand the dynamics of honey bee colonies 
	we present a discrete time marked renewal process model for the size of a colony. We demonstrate that under mild conditions
	 this attains a stationary distribution that depends on the distribution of the numbers of eggs per batch, {the probability an egg hatches} and the distributions of the times between batches and  bee lifetime.
	 This allows an analytic examination of the effect of changing these quantities.
	 We then extend this model to  cyclic annual effects where for example the numbers of eggs per batch and {the probability an egg hatches} may vary over the year.
			\end{abstract}
			
	
	\end{frontmatter}

		\section{Introduction}

The world population reached 7.3 billion as of mid-2015, indicating that the world has added approximately one billion people in the last twelve years and is projected to pass 9 billion by 2050 \citep{WorldPop2015}. It is vital, therefore, that optimal use is made of the world's finite resources for food production \citep{Collier2009}.  Whilst honey bees are not the only pollinators \citep{Rader2015}, it is estimated that 70\% of important food crops are bee-pollinated \citep{Datta2013} and so the status of pollinators, and bees in particular, is of a key concern in global food security \citep{Lonsdorf2009,Blaauw2014}. Many large-scale agricultural businesses now employ peripatetic bee services in order to ensure adequate pollination of crops \citep{Gordon2014, Bishop2016}. There is evidence that pollinator populations are declining and this is a cause for concern for policy-makers \citep{Hafi2012, Breeze2012, Vanbergen2014, EPILOBEE2014, Potts2010}.
This observation motivates our development of a probabilistic model for the colony size. The model allows us to examine how the size depends on the distribution of the numbers of eggs per batch, the chance of hatching and the distributions of the times between batches and  bee lifetime. In turn this gives insight into how a disease for example that affects one or more of these quantities affects the colony size.
		
The dynamics of a bee colony are complex. After her mating flight, the queen remains in the hive laying eggs, fed and cared for by worker bees.  The number of eggs laid are influenced by nectar flow into the colony, which in turn, is influenced by the number of bees out foraging and seasonal effects on availability. During winter, in the absence of forage, there is little or no egg-laying and the colony relies on the honey and pollen stores accumulated during the summer.  {In this situation t}he lifetime of a bee can be significantly longer so that the colony can survive winter \citep{Mattila2001}.
Eggs hatch into larvae which pupate after a few days of feeding \citep{Collins2004} and emerge as juvenile adult bees about three weeks after the egg is laid and living a further 30-50 days. Juvenile adult worker bees mature as they work within the hive until they become foragers after about three weeks after emergence. Colony sizes can range from 20,000 to 100,000 \citep{Goodman2015, Owen2015}. 
 Pressures on a colony include short and long term weather effects, for example,
  the impact of the season on the hatching rate  through changes in relative humidity  was observed in \citep{Al2014}.
 Other effects include
 diseases \citep{Arundel2011, Gordon2014, Bull2012, Furst2014}, pests \citep{Chandler2000, Cook2007, Ryabov2014}, predators \citep{Capri2013}, pesticide, fungicide and herbicide use \citep{Botias2015, Dively2015, Godfray2014, Pettis2013}.
 When these pressures reduce the number of mature foragers, juvenile adults become precocious foragers, that are less efficient and shorter-lived than mature foragers \citep{Schulz1998}.  This can lead to collapse of the colony \citep{PerryMyerscough2015}.  
 We believe that a necessary first step in understanding the  dynamics of bee colonies is to understand their behaviour in a stable environment.

{The literature contains both probabilistic models examined by simulations and deterministic mathematical models that can be used to} examine the effects of environmental changes and disease on colony dynamics \citep{Becher2013, Khoury2013, Lever2014}. These models are themselves complex and, whilst the effects are evident in the results, they often give little understanding of the mechanism. For example, the model description for the simulation program BEEHAVE is comprehensive \citep{Beehave2014} but it is complex
and mathematically intractable.
Our contribution is a parsimonious {and analytically tractable probabilistic} model for the daily colony size which captures the key features and behaviour of the natural system under stable conditions.  This is the first step in the development of more complex models but is of interest as we demonstrate the existence of a stationary distribution for the colony size. In particular, the mean of the stationary distribution is simply expressed in terms of the mean of the batch size, the probability that eggs hatch into bees and the means of bee lifetime and time between batches.
Our model can also be viewed as a type of batch immigration-death model: 
bees arrive in batches of eggs laid by the queen and subsequently die. 
		 
The birth and death process with immigration is well-known in the random process  literature (\citeauthor{GS2001}, \citeyear{GS2001}, pp.~276--278; \citeauthor{Lawler06}, \citeyear{Lawler06}, p.~76).
Classical immigration-death process models are first order Markovian. That is, given the population size at time $t$ the population size at time $t+1$ is independent of the past of the process. 
In our case the daily death rate of the bees depends on the age structure of the colony and is not first order Markovian in this sense. In theory  if there is an upper bound on bee lifetimes one could obtain a Markov property by extending the state space but we consider a more direct approach based on marked renewal processes.

Our approach requires modelling the distributions of four quantities. The first  is the time in days between the days on which the queen lays batches of eggs. Often this will be one day but our results allow this to have a distribution over the positive integers, with daily layings being a special case. 
The second  is  the number of eggs in a batch. The third  is the probability an egg successfully hatches into a bee. The fourth is the lifetime of a bee. This initial model represents a perfect bee world without seasonal effects where there are always sufficient resources so that the batch sizes on different days have the same distribution and the probability that an egg hatches into a bee remains the same. We {then extend this to} a cyclic model where the mean batch size  and the probability of an egg hatching vary over an annual cycle which allows us to consider seasonal effects. 
		 
Our initial hypothesis was that for the given distributions of the time between batches, the number of eggs in a batch, the indicator whether an egg successfully hatches and the lifetime of a bee, with these quantities being mutually independent, there is  a stationary distribution for the size of the colony.

The colony size is modelled as  a marked renewal process with the number of batches being a renewal process and the number of bees hatching in each batch being the mark. Then coupling techniques \cite[pp.~21-34]{Lindvall92} can be used to study its stationary distribution. 
In Section \ref{sec-model} we describe the model and in Section \ref{sec-results-discrete} we give our  theoretical results.
The main result, Theorem~\ref{discretethm1}, gives an explicit representation of the stationary distribution of the {colony size} as  the number of days  becomes large.
To our knowledge such a concise construction for a marked renewal process is not noted elsewhere in the literature.
Moreover, when the distribution of the time between the batches is nonrandom, which is a typically case for a bee colony, the representation also enables us to obtain the characteristic function of the stationary distribution of the colony without relying on the renewal technique. See Corollary~\ref{corollary2}.
			 In Theorem~\ref{cyclicthm1} we give results on the  stationary distribution for the cyclic model.
			 The practical implications of our results  are that the colony size will not endlessly increase but will become stable. In effect it gives a baseline for the colony size. Moreover, our results give an explanation for the recovery in the colony size after swarming. If we suppose that the effect of environmental factors or disease is to either reduce the batch size and the hatching probability, increase mortality through reducing bee lifetime by shifting the lifetime distribution to the left, increase the time between batches or a combination of these effects then {our results allow us to} quantify their effect on the stationary distribution. In Section \ref{sec-applic} we examine this 
		for some parametric models. 
		In Section \ref{sec-sims} we give some simulations to verify the theoretical results and explore extinction probabilities as well as the time taken to attain the stationary distribution. The proofs of the results in Section \ref{sec-results-discrete} are given in the appendices.

	\section{Model and Notation}\label{sec-model}
	
	Let $\tau_i$ be the time between batch $i-1$ and batch $i$, $\zeta_i$ be the number of eggs in batch $i$, $I_{ij}\in\{0,1\}$ be the indicator that the $j$th egg in batch $i$  hatches
 and $\eta_{ij}$ be the lifetime of bee $j$ in batch $i$,  $i\in\N$.  Then $S_i=\sum_{j=1}^i \tau_i$ is the time that the $i$th batch of eggs is
laid. The number of batches laid by time $t$ is $N_t=\sum_{i=1}^\infty \bone_{[S_i,\infty)} (t)$, where $\bone_A(.)$ is the usual indicator function of the set $A$. Clearly, the number of bees that hatch in batch $i$ is
 $\xi_i=\sum_{j=1}^{\zeta_i}I_{ij}$.  Without loss of generality, shifting $\{S_i\}$ by a fixed constant if necessary, we assume that the hatch time from an egg to a bee is zero so that the number of bees in the
colony at time $t$ is $M_t=\sum_{i=1}^{N_t}\sum_{j=1}^ {\xi_i} \bone_{[S_i,S_i+\eta_{ij}]}(t)=\sum_{i=1}^{N_t}\sum_{j=1}^{\xi_i}\bone_{\{S_i+\eta_{ij}\ge t\}}$. 
In
this expository work we suppose that $\{\tau_i,i\in\N\}$ are independent and identically distributed (iid) as are $\{I_{ij},i,j\in\N\}$ and $\{\zeta_i,i\in\N\}$. We further suppose that $\{\eta_{ij}:\ i,j\in\N\}$ are iid and  that the sequences $\{\tau_i,i\in\N\}$, $\{\zeta_i,i\in\N\}$, $\{I_{ij},i,j\in\N\}$ and $\{\eta_{ij},i,j\in\N\}$ are independent.

Under our assumptions the times  $\{S_i,i\in\N\}$ at which the batches  are laid form a renewal process. We let $F_\tau(x)$ and $F_\eta(x)$, $x\ge 0$, denote the distribution functions of $\tau_1$ and $\eta_{11}$ respectively and we suppose that $\mean(\tau_1)$, $\V(\tau_1)$, $\mean(\eta_{11})$, $\V(\eta_{11})$ are all finite. Note that our zero time is arbitrary and in general will not be the time at which a batch of eggs was laid. That is, the time after we start observing until the first batch is laid is $\widetilde \tau_1=\tau_1-u$ for some $u$ representing the (random) time after the last batch of eggs was laid {that  our observations begun}.

	We assume that $\tau_1$ takes positive integer values with distribution $f_j=\prob(\tau_1=j),\ j\ge 1$, where $\sum_{j=1}^\infty f_j=1$. If $\prob(\tau_1=0)>0$, then we can redefine the $\zeta_i$'s so that the new renewal process has inter-renewal times taking positive integer values. Let $X_n$ be the time till the next renewal from time $n$, then $\{X_n\}$ is a Markov chain with transition probabilities $p_{ij}:=\prob(X_{n+1}=j|X_n=i)=f_j$ for $i=1$, $j\ge 1$;
	$p_{ij}=1$ for $i\ge 2$, $j=i-1$; and $p_{ij}=0$ for all other cases. Let $d$ be the greatest common divisor
	of $\{j:\ f_j>0\}$. It is well-known that if $d=1$, then $\{X_n\}$ is ergodic with stationary distribution 
	$$\pi_k:=\frac1\mu(1-F(k-1)),\ k\ge 1,$$
	while for $d>1$, $\{X_n\}$ is not ergodic.
	
	For simplicity, we also assume that the {lifetimes }$\eta_{ij}$'s only take non-negative integer values. We use $\equald$ 
	to stand for two random variables being equal 
	in distribution, $\convd$ for the convergence in distribution and $\sim$ {denotes ``is distributed as''}. We let $r=\prob(I_{11}=1)$ be the probability {an egg hatches} {and set} $q_k=\prob(\eta_{11}=k)$, $k\ge 0$.  Let $\tilde\tau_1$ be an integer-valued random variable {with}  distribution $\{\pi_k, k\ge {1}\}$ such that $\tilde \tau_1$ is independent of $\{\tau_i,\xi_i,\eta_{ij}:\ i,j\ge1\}$. Define
	$\ScS_i=\tilde\tau_1+\tau_2+\dots+\tau_i$, $i\ge 1$.

	\section{Results}\label{sec-results-discrete}

To state the main result, we need a metric to quantify the speed of convergence of $\{M_n\}$ to its stationary distribution. We define the total variation distance between two integer valued random variables $Y_1,Y_2$ as
$$d_{TV}(Y_1,Y_2)=\sup_{B\subset \Z}|\prob(Y_1\in B)-\prob(Y_2\in B)|,$$
where $\Z$ is the space of all integers.

	\begin{theorem} \label{discretethm1} With the setup in \S\ref{sec-model}, let 
	$\ScM:=\sum_{i=1}^\infty\sum_{j=1}^{\xi_i}\bone_{\{\ScS_i\le \eta_{ij}+1\}}.$
	Assume $d=1$ and $\mean\zeta_1<\infty$. Then
	\begin{eqnarray*}
	&&d_{TV}(M_n,\ScM)\\
	&&\le\left\{\begin{array}{ll}
	r\mean(\zeta_1)\mean(\max\{\eta_{11}+1 -n,0\})\le O(n^{-1}),&\mbox{\rm if }F_\tau\mbox{\rm\  is degenerate,}\\
	O\left(n^{-1/2}\right),&\mbox{\rm if }F_\tau\mbox{\rm\  is non-degenerate.}\end{array}\right.\end{eqnarray*}

	Moreover, $\mean \ScM=r\mean(\zeta_1)(\mean\eta_{11}+1)/\mean\tau_1$,
	\begin{eqnarray}
	\V(\ScM)&=&\mean \ScM-(\mean\ScM)^2+r^2(\V(\zeta_1)-\mean\zeta_1)[\mean(\eta_{11}\wedge\eta_{12})+1]/\mean\tau_1\nonumber\\
	&&+\frac{(r\mean\zeta_1)^2}{\mean\tau_1}\sum_{v=0}^\infty\left(2\sum_{i=1}^{v}H(i)+v+1\right)q_v^2
	\nonumber\\
	&&+2\frac{(r\mean\zeta_1)^2}{\mean\tau_1}\sum_{0\le v_1<v_2}\left(\sum_{i=1}^{v_1}H(i)+v_1+1+\sum_{i=v_2-v_1}^{v_2}H(i)\right)q_{v_1}q_{v_2},
	 \label{discretethm2}
	\end{eqnarray}
	where $H(n):=\mean N_n$ satisfies the renewal equation
	\begin{equation}H(n)=F(n)+\sum_{k=1}^nH(n-k)f_k. \label{discretethm3}\end{equation}
	
	\end{theorem}

Theorem~\ref{discretethm1} not only tells us that the marked renewal process $\{M_n\}$ stabilises as $n\to\infty$ but also gives an explicit representation of the stationary distribution. It enables us to simulate the stationary distribution, i.e., the distribution of $M_\infty$, without having to simulate the process $\{M_n\}$ for a long time.

It can be observed from Theorem~\ref{discretethm1} that when $\tau_i$'s are not random ({i.e. $F_\tau$ is degenerate)} and the lifetimes $\eta_{ij}$ are bounded by a constant $K$  then $M_n$ reaches stationarity for $n\ge K+1$. For a bee colony it is reasonable to take $K$ to be 80 days. {An implication of this is that after swarming (where the size of the colony is typically halved) in the stable case that we consider here, the colony size will quickly recover to the stationary size. Thus, at least in the stable case,  no other mechanism is required to explain the recovery of the colony size after swarming. }

	{Let $\Pn(\lambda)$ denote} the Poisson distribution with mean $\lambda$.
 	When $\tau_i$'s are not random, we have the following corollary.
			\begin{corollary}\label{corollary2}
			If $F_\tau$ is degenerate so that we can assume $\ScS_i,\ i\ge 1,$ are constants, then the stationary distribution of $\{M_n\}$ has the characteristic function
		\begin{equation}\phi_{\ScM}(u):=\mean e^{\ime u\ScM}=\prod_{i=1}^\infty\psi_\zeta\left(r(e^{\ime u}-1)(1-F_\eta(\ScS_i-2))+1\right),\label{cor2-1}\end{equation}
		where $F_\eta$ is the distribution function of $\eta_{11}$ and $\psi_\zeta(s):=\mean s^{\zeta_1}$. In particular, if $\zeta_1\sim \Pn(\lambda)$, then $\ScM\sim \Pn\left(r\lambda\sum_{i=1}^\infty\left(1-F_\eta(\ScS_i-2)\right)\right).$
		\end{corollary}

	{If the minimum size of a viable colony is $\nu$ {then
	using Corollary~\ref{corollary2}} to approximate $\prob(M_n<\nu)$ gives an approximate  lower bound on the probability of colony extinction at time $n$. This is a lower bound as the colony may have become extinct before $n$ (i.e. $\prob(\min_{k\le n} M_{{k}} <\nu) \ge \prob(M_n<\nu)$). The computation of the extinction probabilities requires derivation of the distribution of
	$\min_{k\le n} M_{{k}}$ which is beyond the present work and for now we address this through simulations in Section \ref{sec-sims}. See Figure \ref{fig:extinction_prob}.}

	To {allow} the distribution $F_{\zeta,i}$ of $\zeta_i$, the hatch probability $r_i$ of eggs in batch $i$ and the distribution $F_{\eta,i}$ of $\eta_{ij},\ j\ge 1$, {to all} depend on time $i$, 
	we propose the periodic model $\{M_{lD+i}, i=1,\dots,D \}$, $l\in\N$, to mimic the behaviour over a {cycle} of $D$ days, {where} $D=365$ days represents a year. That is we assume $F_{\zeta,\ell}=F_{\zeta,i}$, $r_{\ell}=r_i$ and $F_{\eta,\ell}=F_{\eta,i}$ if $\ell=i\mod D$.

 		\begin{theorem} \label{cyclicthm1}
	Assume $\tau_1\equiv 1$, $\mean \zeta_i<\infty$ for all $1\le i\le D$, and $\prob(\eta_{ij} \leq K) = 1$, then $(M_{lD+n})_{l\in\mathbb{N}}$ reaches stationarity for $lD+n>K$, and the characteristic function of the stationary distribution is
	$$\phi_{n+lD}(t) = \prod_{i=n+lD-K}^{n+lD}\psi_{\zeta,i}\left(r_i(e^{\ime t}-1)(1-F_{\eta,i}(n+lD-i-1))+1\right),$$
	where $\psi_{\zeta,i}(s):=\mean s^{\zeta_i}$.
	In particular, $\mean M_{n+lD} = \sum_{i=n+lD-K}^{n+lD}(1-F_{\eta,i}(n+lD-i-1))r_i\mean \zeta_i$, and if $\zeta_i$ follows $\Pn(\lambda_i)$, then\\ $M_{n+lD} \sim  \Pn\left(\sum_{i=n+lD-K}^{n+lD}r_i\lambda_i(1-F_{\eta,i}(n+lD-i-1))\right)$.
	\end{theorem}	
	
This theorem establishes that the year to year behaviour of the colony size on a given day is stationary.

	\begin{figure}[t]
	  \centering
	    \includegraphics[width=.7\textwidth]{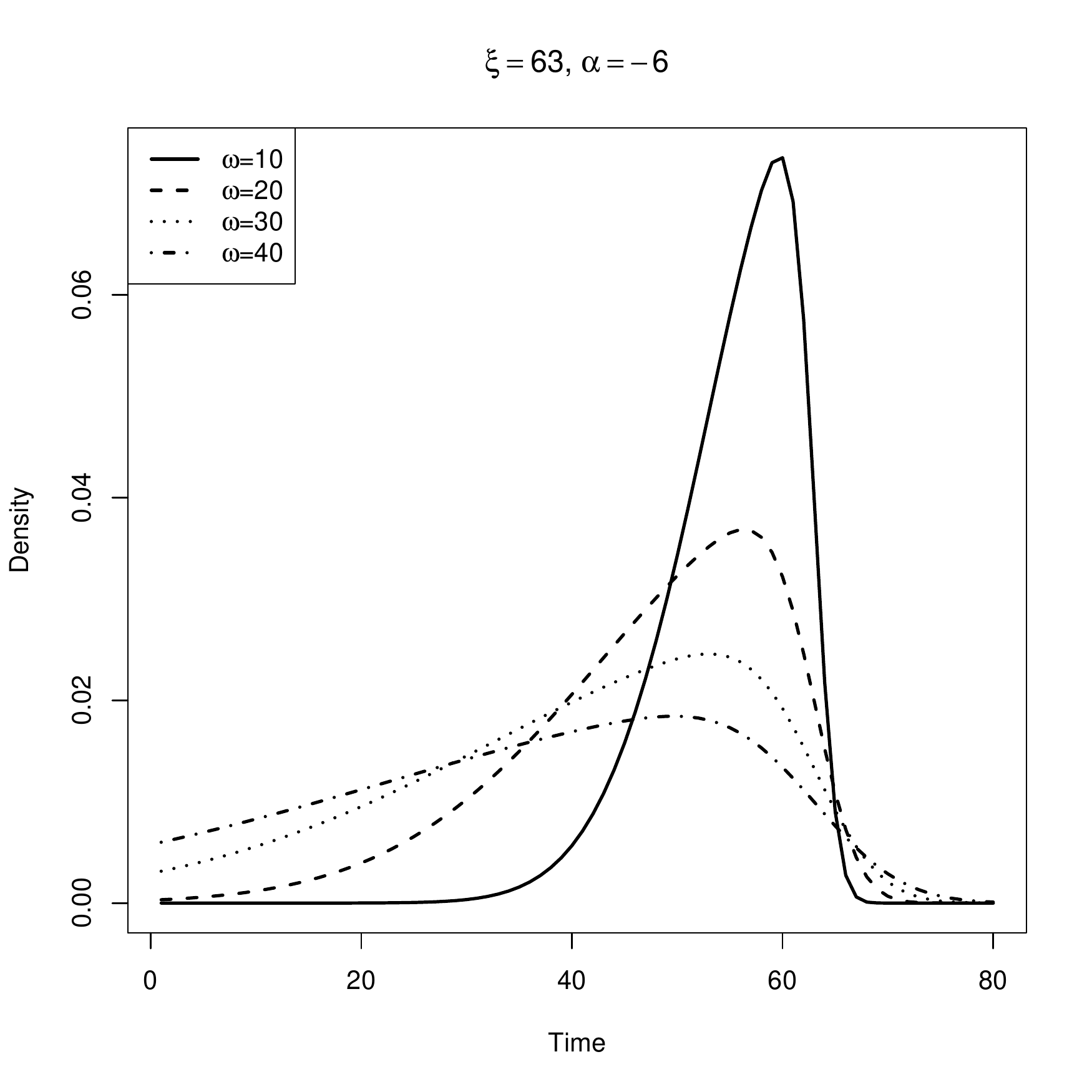}
	  \caption{The $SN(63,\omega^2,-6)$ density for  $\omega \in \{ 10,20,30,40\}$.}
	  \label{fig:sn}
	\end{figure}

	\section{Applications}\label{sec-applic}

	We first use our analytic results to examine the effects of changing  the distributions of  the hatch probability, the  batch size, the bee lifetime and the times between batches on the stationary distribution of the colony size. 
	For example with a mean batch size of 1875,  hatching probability $0.8$,  mean bee lifetime of 63 days and a constant time between batches of 1, then from Theorem \ref{discretethm1} the stationary mean colony size is 96,000. If the mean batch size drops to 1600 and the hatch probability drops to $0.75$ while the mean bee lifetime remains at 63 days,  the stationary mean colony size is now 76,800. If the mean bee lifetime also drops to 50 , the stationary mean colony size is 61,200.
	Thus if the effect of a factor is to reduce the nectar flow into the hive reducing the mean batch size our results allow us to quantify this effect on the mean colony size. A reduction in the nectar flow could occur from a change in land use, change in climate, or a reduction in the numbers of forager bees.
	Alternatively, a disease may directly reduce the batch size or  the hatching probability of the eggs.

	To examine the behaviour of the colony {size in more detail}
    we model the lifetimes using a discretization of the skew normal distribution \citep{Azzalini2013}  truncated at zero as this distribution allows left skewed lifetime distributions.
	Left skewness is important as mortality will be low until bees start foraging and this may not occur until they are aged 40 days or so.
   	    The skew normal is a three parameter distribution with a location parameter $\Xi$,
	     a scale parameter $\omega$ and a slant parameter $\alpha$ \citep{Azzalini2013}. Its density is
   	  \[
   	  f(x;\Xi,\omega,\alpha)=\frac{2}{\omega}\phi\left(\frac{x-\Xi}{\omega}\right)\Phi\left(\alpha \frac{x-\Xi}{\omega}\right),
   	  \]
	  where $\phi$ and $\Phi$ respectively denotes the probability density function and the cumulative distribution function of the standard normal law.
	   Note that when $\omega = 0$ the distribution of $\eta$ is not random but constant equal to $\Xi$.
   	  We denote this density $SN(\Xi,\omega^2,\alpha)$. The mean of the skew normal distribution is $\mu=\Xi+\omega\delta\sqrt{2/\pi}$ where $\delta=\alpha/\sqrt{1+\alpha^2}$.
   	  If $\alpha<0$ then this distribution is skewed to the left. We discretise this distribution by rounding up to the nearest integer, with negative values rounded up to one.
	  
	  		  	\begin{figure}[t]
	  	  \centering
	  	    \includegraphics[width=.9\textwidth]{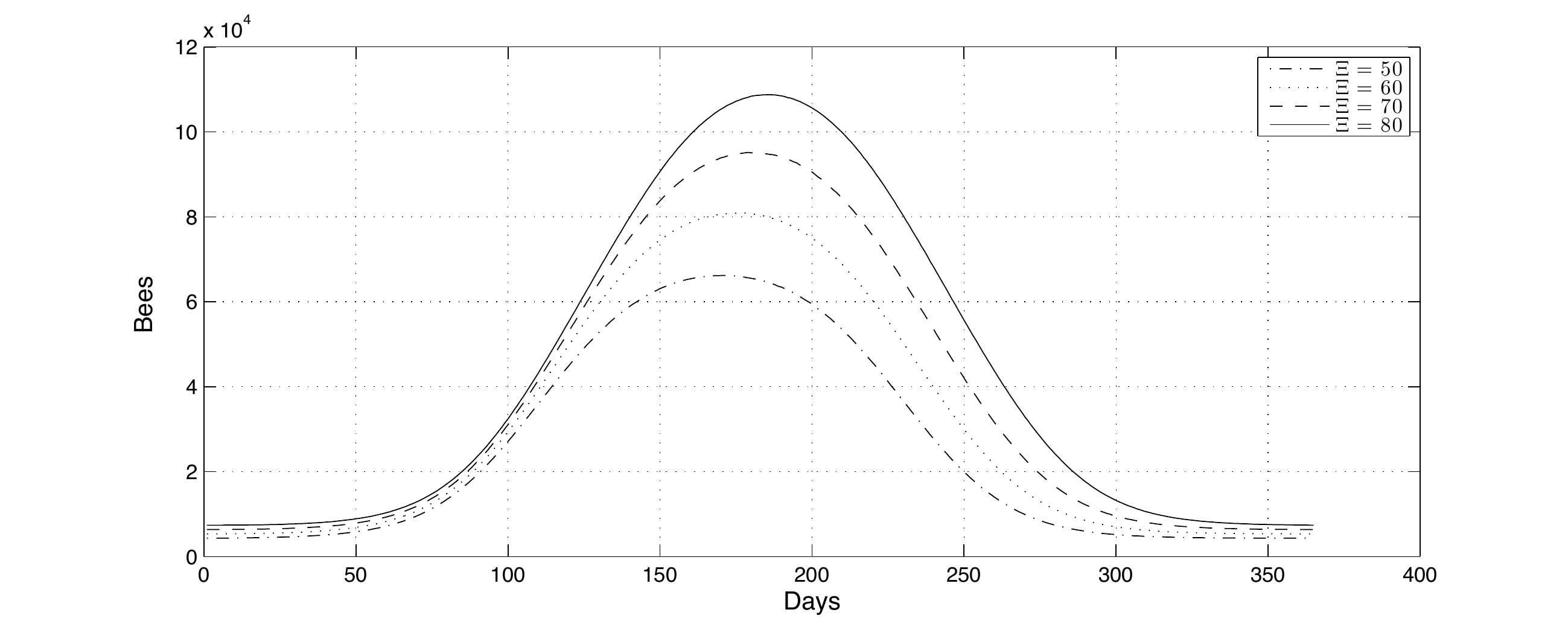}
	  	  \caption{The mean of the daily colony size $M_n$ when the lifetime $\eta \sim SN(\Xi,10^2,-6)$, for the location parameter $\Xi \in \{50,60,70,80\}$ and a cyclic mean number of eggs
		  given by $\mean(\xi_i) {=r_i\mean(\zeta_i)}= 1500\times seas(i) + 100$.}
	  	  \label{fig:cyclic-mean-xi}
	  	\end{figure}
	  
	  In Figure \ref{fig:sn} as an illustration we plot the $SN(63,\omega^2,-6)$ density for varying values of $\omega$. This indicates that the shape of the distribution can become quite skewed.
	  This skewness gives one way of representing an increase in bee mortality, with an increasingly heavy left tail when more bees die at a younger age. Varying the location and slant parameters further change the shape and location of this distribution, which gives the model considerable flexibility. Therefore we use this model to illustrate our results.

{For convenience, we assume that we start observations on the first day of the year so that the subindex $i$ represents the $i$th batch as well as the $i$th day (replace $i$ with $i\mod365$ if $i>365$) of the year.	
  	
In Figure~\ref{fig:cyclic-mean-xi} we use Theorem \ref{cyclicthm1} and the skew-normal model for the lifetimes  to investigate the effect of increasing the location parameter $\Xi$ 
 when there is a cyclic
model for numbers of eggs laid. 
That is, $\mean(\xi_i) = r_i\mean(\zeta_i)=1500 \times seas(i) + 100$ 
where $seas$ is given by the formula of \citet[p.~3]{Schmickl2007}.
We suppose a Poisson distribution for $\zeta_i$ and have noted that by the thinning property of the Poisson distribution
(\citeauthor{GS2001}, \citeyear{GS2001}, p.~255 and p.~287), $\xi_i$ also follows a Poisson distribution with mean $\mean(\xi_i)=r_i\mean(\zeta_i)$. 
This shows that increasing the location parameter $\Xi$ increases the colony size.

In  Figure  \ref{fig:cyclic-mean-w} in \ref{app-plots} we again use the skew-normal model for the lifetimes and
in these {figures} we examine the effect of varying the other parameters in the lifetime distribution on the means as well as the seasonal effects.
In Figure   \ref{fig:cyclic-mean-w} a)  we see that increasing $\omega$ decreases the mean seasonal population size. This is consistent with Figure \ref{fig:sn}
where we saw that increasing $\omega$ resulted in a heavier left tail to the lifetime distribution.
Figure  \ref{fig:cyclic-mean-w} b)
illustrates the effect of changing values of $\alpha$ and  
 in Figure  \ref{fig:cyclic-mean-w} c), we plot the mean number of bees over a year with several distributions for $\zeta_i$ and varying hatch probabilities $r_i$.
 As one would expect the mean population size increased as the mean number of eggs laid increases.

\section{Simulations}\label{sec-sims}

  		  	\begin{figure}[t]
  	  \centering
  	    \includegraphics[width=.9\textwidth]{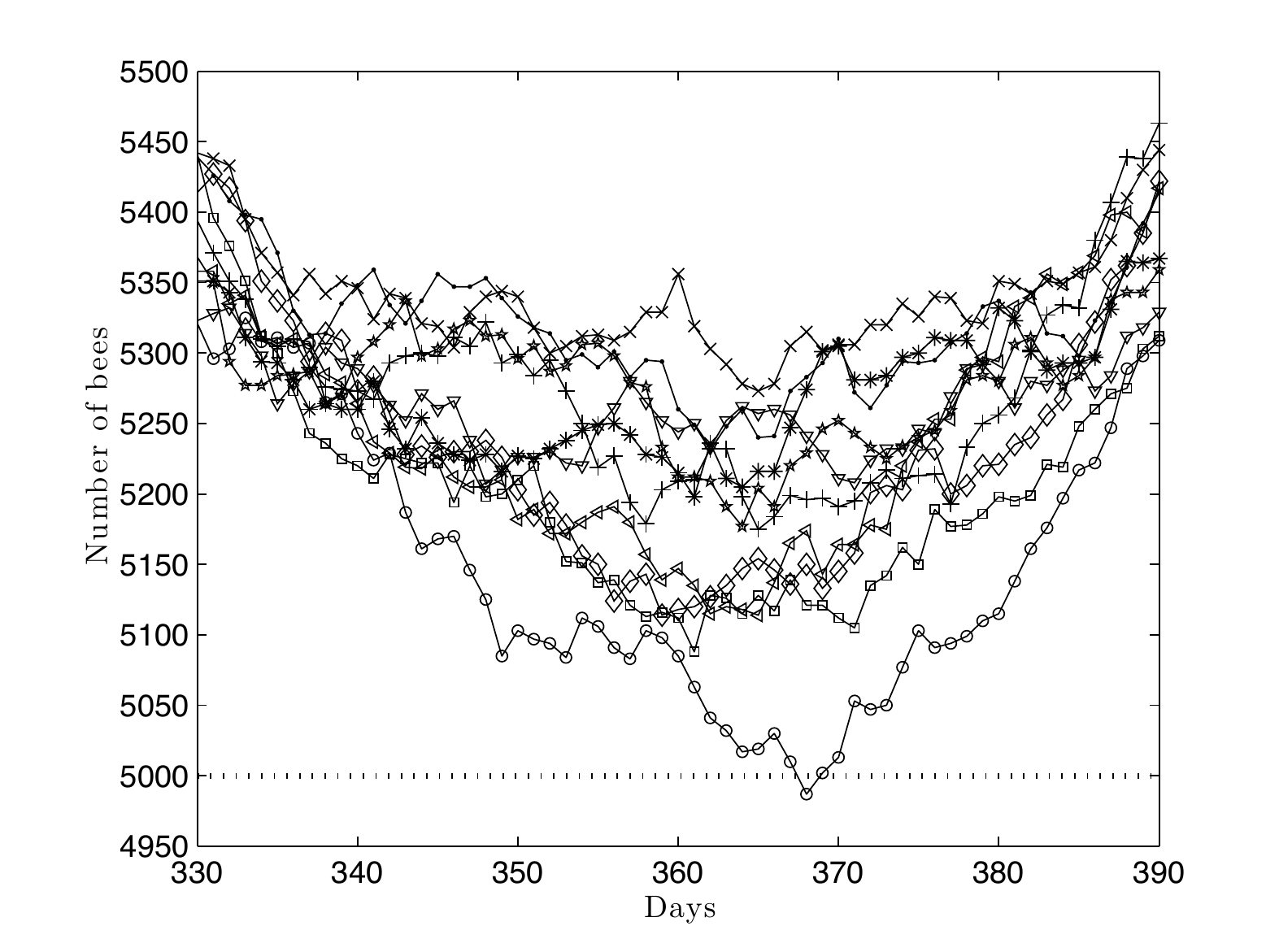}
  	  \caption{The number of bees during the first winter over ten different simulation, when the lifetime $\eta \sim SN(63,30^2,-6), $ and \{the number of hatched bees $\xi_i\sim\Pn(1500 \times seas(i) + 127)$. The dotted line is the threshold of $5000$ under which the colony collapses.}
  	  \label{fig:simulation_winter}
  	\end{figure}
	
	To verify the theoretical results we conduct a series of simulations. We initially {suppose} $\zeta_i \sim \Pn(\mu_\zeta)$ and as the  probability {of an egg hatching} is $r$ we have $\xi_i\sim \Pn(\mu_\xi)$ with $\mu_\xi=r\mu_\zeta$.
In Figure~\ref{fig:simulation_winter} we report ten simulations of our model with $\eta \sim SN(63,30^2,-6), $ and $\mean(\xi_i) =r_i\mean(\zeta_i) = 1500 \times seas(i) + 127$ during the first winter. We can see that one of the simulations resulted in the number of bees being under the nominal threshold of $5000$ and in that case the colony would have become extinct.
To examine this further if the laying times are constants and the {distribution} of the numbers of eggs laid, the hatching probability and the distribution of the bee lifetimes are specified, we can apply Corollary~\ref{corollary2} to determine lower bounds on the extinction probabilities. In general, we can always use simulations to estimate extinction probabilities, although for many parameter values this will be small. In Figure~\ref{fig:extinction_prob} we plot the extinction probability (supposing that extinction occur when the number of bees in a colony is less than $5000$) according to the mean number of eggs laid per day.
{Here, we simulated the probability of extinction over 30 days as a function of the mean number of eggs laid. This number corresponds to mid-winter in our model.}  {Thus our model may be used to determine the effects of changes in the underlying parameters on the extinction probabilities.} Note that it only requires a small change in the mean number of eggs laid per day for the extinction probability to increase from zero to one.

	  		  	\begin{figure}[t]
  	  \centering
  	    \includegraphics[width=.9\textwidth]{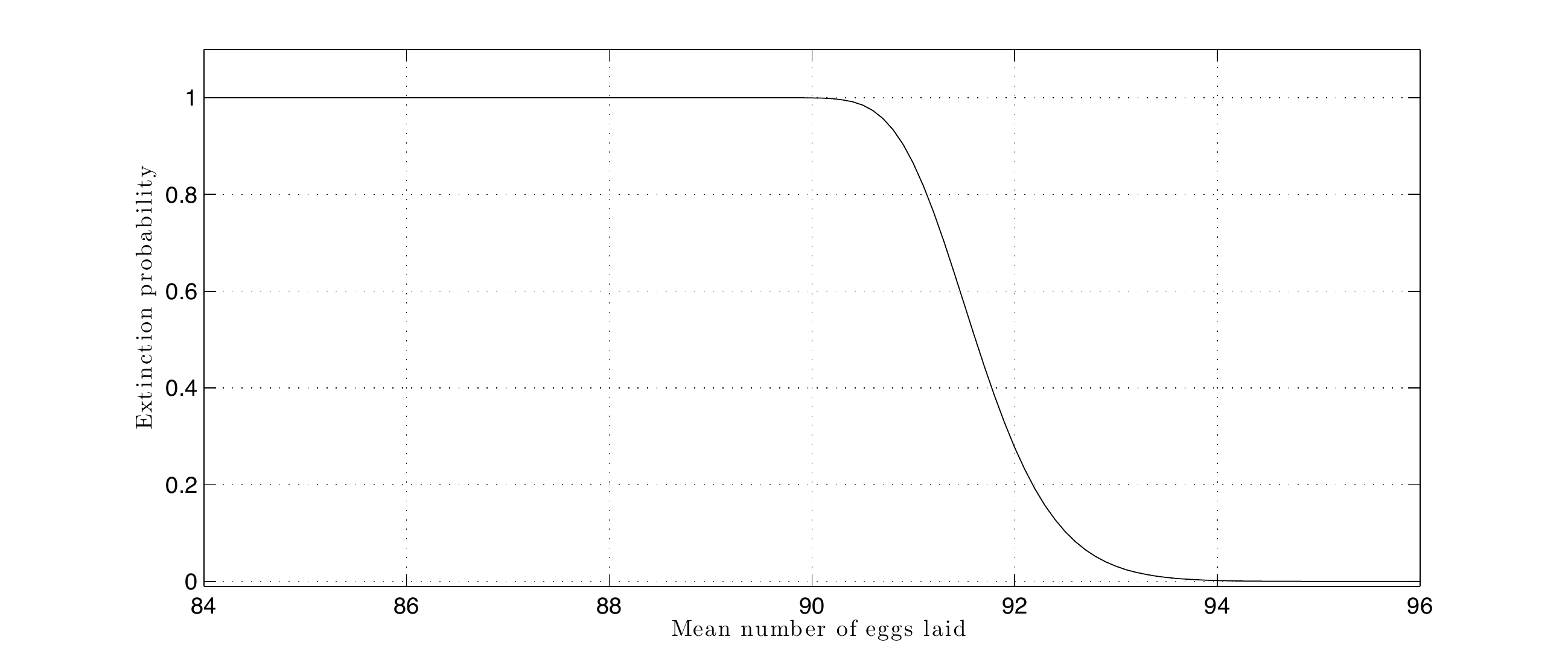}
  	  \caption{The {simulated} extinction probability as a function of the mean number of eggs laid per day obtained in Corollary \ref{corollary2}. Here we took $ \eta \sim SN(63,10^2,-6)$ and $\xi_i\sim \Pn(L)$ with $L\in[80;95]$.}
  	  \label{fig:extinction_prob}
  	\end{figure}
	   
In Figure~\ref{fig:swarming} we plot the number of bees in the colony during a {period slightly longer} than one year {using the model described above but incorporating swarming. That is} we assumed that when the {colony size}  was more than $80000$, half of the colony leave the colony. Clearly, in the right conditions the colony size recovers {and as noted above no other mechanism is required to explain the recovery of the colony size. }

	   A further quantity of interest is the time until the colony size attains the stationary distribution. In our case, as the maximum lifetime of a bee is bounded by $80$, the colony size attains the stationary distribution after {no more than $81$} days as long as the number of eggs laid is sufficiently large. This can be observed in Figure \ref{fig:swarming}. Indeed, $80$ days after the colony splits, the number of bees in the hive is as if swarming never happened. From an evolutionary perspective by splitting into two, colonies that swarm have an increased probability that one the colonies will survive.
	   
   		  		  	\begin{figure}[ht]
     	  \centering
     	    \includegraphics[width=.9\textwidth]{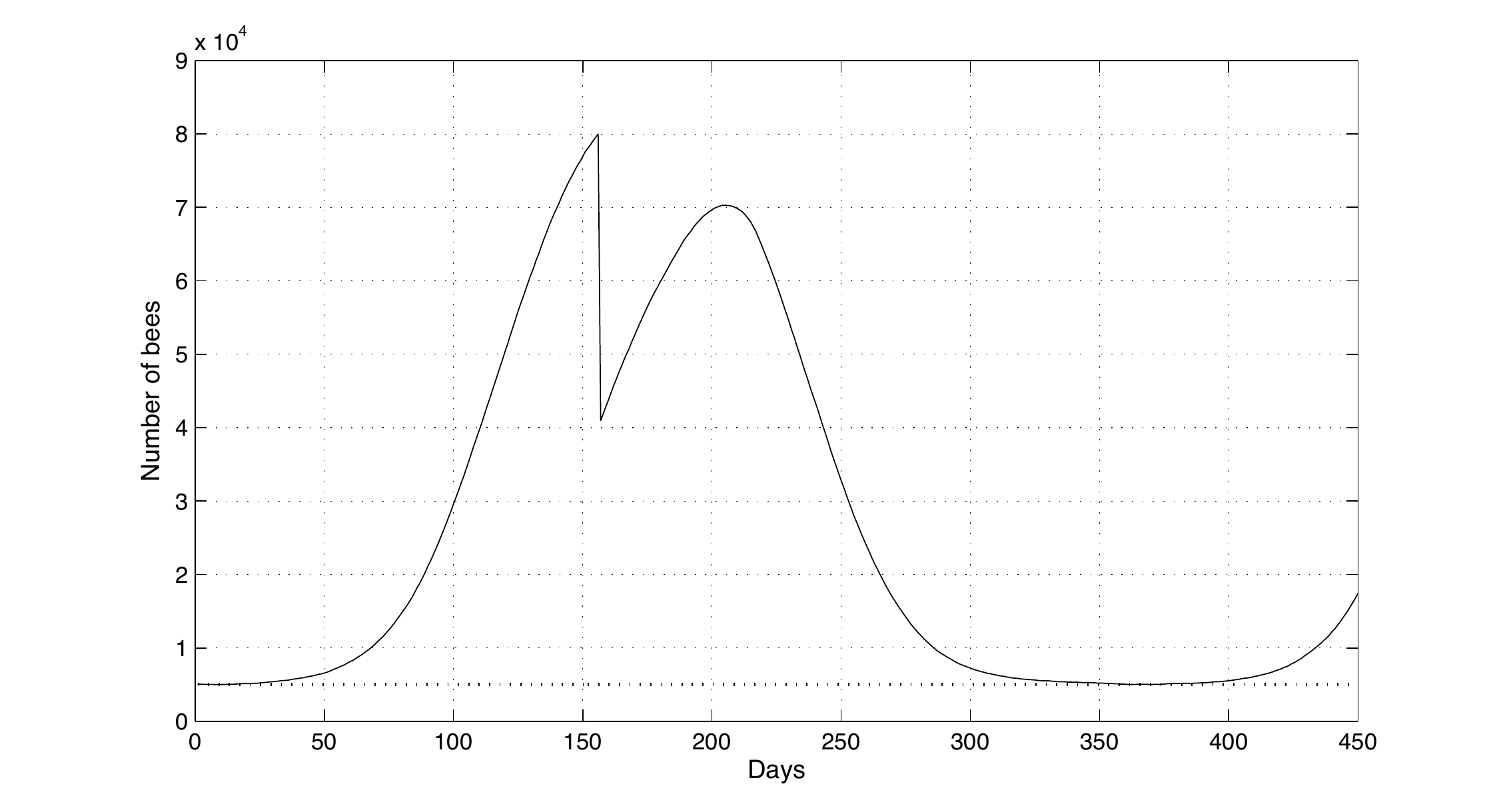}
     	  \caption{The number of bees simulated over 450 days according to the cyclic model with swarming. Here $ \eta \sim SN(63,10^2,-6),$ and {$\xi_i\sim\Pn(1500 \times seas(i) + 90)$} and  we allowed swarming when the number of bees exceeded 80000.}
     	  \label{fig:swarming}
     	\end{figure}

	\section{Discussion}
	
	{Our model is an encouraging  first step in modelling the dynamics of a bee colony. We have developed a probabilistic model for the dynamics of a bee colony in an idealised bee world, which may not of course be attained. 
	This degree of simplicity was necessary to obtain analytic results, but
	nevertheless, our results explain why a typical colony size in summer is around 90,000 bees. Our interpretation is that the stationary distribution represents a state attainable under ideal conditions and thus in many ways is an upper bound. {We extended the stationarity results to a more naturalistic seasonal model and showed the marginal distribution of the colony size on a given stage of the season was stationary. Our model also allowed us to give a lower bound for extinction probabilities and  demonstrate how the colony size recovers after swarming, as seen in the natural system. 
	 The simple expression for the mean colony size is easily computed in terms of the hatching probability and the mean of the numbers of eggs laid, the expected bee lifetime and the expected time between the laying of batches of eggs.
	 A determination of how a {factor} such as environmental perturbations, disease pressure or weather impacts the colony size (for example through increased mortality or reduced nectar flow and hence {the number of eggs laid per day}) then allows an examination of its effect on colony size.

	There are  several ways in which our work can be extended. {An extension of our results would first be to consider the joint behaviour of the colony size over several time points which would potentially allow the derivation of exact extinction probabilities in addition to the lower bound. That is to consider the colony size $\{M_{l D+i},i=1,\dots,D\}$ over a year as ``functional data'' \citep{Ramsay97functionaldata} and examine the behaviour of this function. Corollary~\ref{corollary2} examines the marginal distributions of this quantity. One could then consider the effects of long term trends on the colony size and survival. In particular conditional probabilities of survival given a particular colony size could in theory be derived. That is, given the current colony size of say 7000 bees at a given time of year under what conditions would the probability of surviving over winter be 90\% for example? Such calculations could indicate whether supplementary feeding was required to increase the  {number of eggs laid per day} to help ensure survival. By relating the effects of environmental changes and disease to these quantities, their effect on colony size may also be examined.
	Whilst our current results are for stable and cyclic effects we anticipate extending these to even more complex and realistic situations in the future. For example, the availability of nectar in the hive may depend on the numbers of forager bees so that the average batch size could be taken to be dependent on the number of forager bees. In their deterministic models, \citet{Khoury2013}
	give the steady state population as a function of bee death rate and food collection rate. Bee death rate relates to our bee lifetime distribution. However, we can consider more subtle effects. We do not explicitly include food collection rate in the current model but we have argued this could be included in the models for the mean batch size.
	The availability of nectar could also depend on the age related efficiency of the forager bees and the average batch size could be taken to be dependent on the ages of the forager bees. These more complex feedback models are technically challenging and beyond the scope of the current work.
	
	Regarding more technical extensions, recall that a lattice distribution is the distribution of a random variable taking possible values of the form $\{a+bk, a,b\mbox{ are constants}, b\ne0, k\in\Z\}$. The distribution $F_\tau$ of the renewal times $\{\tau_i\}$ in Theorem~\ref{discretethm1} is a special case of a lattice distribution with $b=1$.
	One can easily adapt the proof of Theorem~\ref{discretethm1} to obtain the stationary distribution when $F_\tau$ either follows a general lattice distribution or a non-lattice distribution.
	We also note that one can also study $\{M_n\}$ as a discrete-time queueing system with batch arrivals (see \cite{Takagi93}).

	We conclude that any change in the parameters of the lifetime distribution or the numbers of eggs laid can affect the population size. This allows a variety of ways in which a disease for example can affect the lifetime distribution and hence colony size. Whilst their properties can be  difficult to derive, we believe that probabilistic models for the size of a bee colony are mathematically tractable and provide useful and readily interpretable information on the dependence of the colony size on the times between batches of eggs being laid, the sizes of the batches, the hatching probability and bee lifetime. We further believe they and their extensions have considerable potential both in understanding the behaviour of bee colonies and also provide new classes of stochastic models worthy of further research.

	\bibliographystyle{chicago} 
	\bibliography{Bees-1}
	
	\clearpage

\begin{appendix}
	
	\section{Proof of Theorem 1}
	
	
	First of all, since $\xi_1=\sum_{j=1}^{\zeta_1}I_{1j}$ is a mixed binomial, we have
	\begin{equation}\mean\xi_1=\mean\left\{\mean\left(\left.\sum_{j=1}^{\zeta_1}I_{1j}\right|\zeta_1\right)\right\}=\mean(r\zeta_1)=r\mean\zeta_1\label{convrate05}\end{equation}
	and
	\begin{equation}\V(\xi_1)=\mean\left\{\V\left(\left.\sum_{j=1}^{\zeta_1}I_{1j}\right|\zeta_1\right)\right\}+\V\left\{\mean\left(\left.\sum_{j=1}^{\zeta_1}I_{1j}\right|\zeta_1\right)\right\}=r(1-r)\mean\zeta_1+r^2\V(\zeta_1).\label{convrate06}\end{equation}

	If $F$ is non-degenerate, we write $T_n=\tau_2+\dots+\tau_n$, then the condition $d=1$ implies $d_{TV}(T_n,T_n+1)=O\left(n^{-1/2}\right)\to 0$ as $n\to\infty$ \cite[pp.~41--43]{Lindvall92}. This ensures
	\begin{eqnarray}&&d_{TV}(\ScS_n,S_n)\le d_{TV}(\ScS_n,T_n)+d_{TV}(T_n,S_n)\nonumber\\
	&&\le (\mean\tilde\tau_1+\mean\tau_1)d_{TV}(T_n,T_n+1)=O\left(n^{-1/2}\right)\to 0\label{discretethm4}
	\end{eqnarray}
	as $n\to\infty$. If $F_\tau$ is degenerate, then $\ScS_n=S_n$, giving $d_{TV}(\ScS_n,S_n)=0$.
	Hence, we obtain
	\begin{equation}d_{TV}(\ScS_n,S_n)=\left\{\begin{array}{ll}
	0,&\mbox{\rm if }F_\tau\mbox{\rm\  is degenerate,}\\
	O\left(n^{-1/2}\right),&\mbox{\rm if }F_\tau\mbox{\rm\  is non-degenerate.}\end{array}\right.\label{convrate00}
	\end{equation}

	Let $\ScN_n:=\sum_{i=1}^\infty\bone_{\{\ScS_i\le n\}}$ and 
	$\ScM_n:=\sum_{i=1}^{\ScN_n}\sum_{j=1}^{\xi_i}\bone_{\{\ScS_i+\eta_{ij}\ge n\}}$, we first establish that
	\begin{equation}
	d_{TV}(M_n,\ScM_n)=\left\{\begin{array}{ll}
	0,&\mbox{\rm if }F_\tau\mbox{\rm\  is degenerate,}\\
	O\left(n^{-1/2}\right),&\mbox{\rm if }F_\tau\mbox{\rm\  is non-degenerate.}\end{array}\right.\label{convrate01}
	\end{equation}
	In fact, when $F_\tau$ is degenerate, we have $M_n\equald\ScM_n$, so \Ref{convrate01} is obvious. For the case that $F_\tau$ is non-degenerate, we have $\sigma_\tau:=\sqrt{\V(\tau_1)}\in (0,\infty)$. Take $a_n=\lfloor 0.75n\rfloor$, $b_n=\lfloor 0.5n/\mean\tau_1\rfloor$, where $\lfloor x\rfloor$ is the biggest integer less than or equal to $x$. By \cite[pp.~253--254]{BHJ92}, enlarging the probability space if necessary, we can construct iid $\{\tau_i'\}$ having the same distribution as that of $\tau_1$ and independent of $\{\xi_i,\eta_{ij}:\ i,j\in\N\}$ such that $\left\{S_k':=\sum_{j=1}^k\tau_j':\ k\in\N\right\}$ satisfies
	\begin{equation}\prob(S_{b_n}'\ne \ScS_{b_n})=d_{TV}(S_{b_n},\ScS_{b_n})=O(b_n^{-1/2})=O(n^{-1/2}),
	\label{coupling01}
	\end{equation}
	where the second equality is due to \Ref{convrate00}. Set $\ScS_j':=\left\{\begin{array}{ll}\ScS_j,&\mbox{\rm if }j\le b_n,\\
	\ScS_{b_n}+\sum_{l=b_n+1}^j\tau_l',&\mbox{\rm if }j>b_n,
	\end{array}\right.$ $N_n':=\sum_{i=1}^\infty \bone_{\{S_i'\le n\}}$, $M_n':=\sum_{i=1}^{N_n'}\sum_{j=1}^{\xi_i}\bone_{\{S_i'+\eta_{ij}\ge n\}}$; $\ScN_n':=\sum_{i=1}^\infty \bone_{\{\ScS_i'\le n\}}$, $\ScM_n':=\sum_{i=1}^{\ScN_n'}\sum_{j=1}^{\xi_i}\bone_{\{\ScS_i'+\eta_{ij}\ge n\}}$, then $M_n'\equald M_n$, $\ScM_n'\equald\ScM_n$ and
	$$d_{TV}(M_n,\ScM_n)\le \prob\left(M_n'\ne\ScM_n'\right).$$
	Therefore, it suffices to show that $\prob\left(M_n'\ne\ScM_n'\right)=O(n^{-1/2})$.
	To this end, we observe that
	\begin{equation}
	\{M_n'\ne \ScM_n'\}
	\subset\{S_{b_n}'\ne \ScS_{b_n}\}\cup\{S_{b_n}'>a_n\}\cup\{S_{b_n}'= \ScS_{b_n}\le a_n,M_n'\ne \ScM_n'\}.\label{coupling02}
	\end{equation}
	For the event $\{S_{b_n}'\ne \ScS_{b_n}\}$, \Ref{coupling01} ensures that $\prob(S_{b_n}'\ne \ScS_{b_n})=O(n^{-1/2})$. For the event $\{S_{b_n}'>a_n\}$, we have	$$
	\prob(S_{b_n}'>a_n)=\prob\left(\frac{S_{b_n}'-\mean S_{b_n}'}{\sqrt{\V(S_{b_n}')}}>c_n\sqrt{n}\right)=O(n^{-1/2})+\int_{c_n\sqrt{n}}^\infty \frac{e^{-x^2/2}}{\sqrt{2\pi}}dx=O(n^{-1/2}),
	$$
	where $c_n\asymp \frac{\sqrt{2\mean\tau_1}}{4\sigma_\tau}$, the second equality is due to \cite[Theorem~11.2]{CGS11} and the last equality follows from \cite[(2.82)]{CGS11}.
	Therefore, it remains to show that the final event of \Ref{coupling02} has probability of order no more than $O(n^{-1/2})$. This can be worked out as follows.
	On $\{S_{b_n}'= \ScS_{b_n}\le a_n,M_n'\ne \ScM_n'\}$, the two marked renewal processes coincide from time $a_n$, so $M_n'\ne \ScM_n'$ can only happen when at least one of the bees before $a_n$ are still alive at time $n$. That is, the bee must have lived through the interval $[a_n,n]$ with lifetime $\eta_{ij}\ge n-a_n$,  giving
	$$\{S_{b_n}'= \ScS_{b_n}\le a_n,M_n'\ne \ScM_n'\}\subset \cup_{i=1}^{b_n-1}\cup_{j=1}^{\xi_n}\{\eta_{ij}\ge n-a_n\}.$$
	This yields
	\begin{eqnarray*}
	&&\prob(S_{b_n}'= \ScS_{b_n}\le a_n,M_n'\ne \ScM_n')\\
	&\le&(b_n-1)\mean\xi_1\prob(\eta_{11}\ge n-a_n)
	\le b_n\mean\xi_1\int_{n-a_n}^\infty \frac{x^2}{(n-a_n)^2}dF_{\eta_{11}}(x)\\
	&\le&\frac{b_nr\mean\zeta_1\mean\eta_{11}^2}{(n-a_n)^2}=O(n^{-1}),\end{eqnarray*}
	where the last inequality is due to \Ref{convrate05}.
	This completes the proof of \Ref{convrate01}.
	
		Set $\ScX_n:=\ScS_{\ScN_n+1}-n$, $\ScB_n:=n-\ScS_{\ScN_n}$, then $\ScX_n$ is the time to the next renewal from time $n$ and $\ScB_n$ is the time from the renewal at or before time $n$ to time $n$. Clearly $\ScX_0=\tilde\tau_1\sim\{\pi_k\}$. For $l\in\N$,
	\begin{eqnarray*}
	\prob(\ScX_1=l)&=&\prob(\ScX_1=l|\tilde\tau_1=1)\prob(\tilde\tau_1=1)+\prob(\ScX_1=l|\tilde\tau_1=l+1)\prob(\tilde\tau_1=l+1)\\
	&=&f_l\pi_1+\pi_{l+1}=\pi_l,
	\end{eqnarray*}
	hence $\ScX_1\sim\{\pi_k\}$. By the induction on $n$, $\ScX_n$ has the distribution $\{\pi_k\}$ for all $n$. Next, we establish that $\ScB_n+1$ follows the distribution $\{\pi_k\}$ as well. To this end,
	for $l\in\N$,
	$$\prob(\ScB_n\ge l)=\prob(\ScN_n-\ScN_{n-l}=0)=\prob(\ScX_{n-l}\ge l+1)=\sum_{k\ge l+1}\pi_k,$$
	so the claim follows immediately. A direct consequence of the claim is that
	$$n+1-(\ScS_{\ScN_n},\ScS_{\ScN_n-1},\ScS_1)\equald(\ScS_1,\ScS_2,\dots,\ScS_{\ScN_n}),$$ 
	which ensures 
	$$\ScM_n\equald\sum_{i=1}^{\ScN_n}\sum_{j=1}^{\xi_i}\bone_{\{n+1-\ScS_i+\eta_{ij}\ge n\}}
	=\sum_{i=1}^{\ScN_n}\sum_{j=1}^{\xi_i}\bone_{\{1+\eta_{ij}\ge \ScS_i\}}=:\tilde\ScM_n.$$
	Clearly,	
	\begin{equation}
	d_{TV}(\ScM_n,\ScM)\le \prob\left(\tilde\ScM_n\ne \ScM\right)\le \mean(\ScM-\tilde\ScM_n)=\mean\sum_{i=\ScN_n+1}^\infty
	\sum_{j=1}^{\xi_i}\bone_{\{1+\eta_{ij}\ge\ScS_i\}}.\label{convrate02}
	\end{equation}
	However, $\ScS_{\ScN_n+l}\ge n+l$, it follows from \Ref{convrate02} that
	\begin{eqnarray}
	d_{TV}(\ScM_n,\ScM)&\le& \mean\sum_{i=n+1}^\infty
	\sum_{j=1}^{\xi_i}\bone_{\{1+\eta_{ij}\ge i\}}\nonumber\\
	&=&\mean\xi_1\mean(\max\{\eta_{11}+1 -n,0\})	\label{convrate03}\\
	&\le& \mean\xi_1\mean\left(\eta_{11}\bone_{\{\eta_{11}\ge n\}}\right)
	\le n^{-1}\mean\xi_1\mean\left(\eta_{11}^2\right)\nonumber\\
	&=&O\left(n^{-1}\right).
	\label{convrate04}
	\end{eqnarray}
		Combining \Ref{convrate01}, \Ref{convrate03}, \Ref{convrate04} and \Ref {convrate05} gives
	\begin{eqnarray*}
	&&d_{TV}(M_n,\ScM)\le d_{TV}(M_n,\ScM_n)+d_{TV}(\ScM_n,\ScM)\\
	&&\le\left\{\begin{array}{ll}
	r\mean(\zeta_1)\mean(\max\{\eta_{11}+1 -n,0\})\le O(n^{-1}),&\mbox{\rm if }F_\tau\mbox{\rm\  is degenerate,}\\
	O\left(n^{-1/2}\right),&\mbox{\rm if }F_\tau\mbox{\rm\  is non-degenerate.}\end{array}\right.\end{eqnarray*}

	For $\mean\ScM$, we have from \Ref{convrate05} that
	\begin{eqnarray*}
	\mean\ScM&=&\sum_{i=1}^\infty\mean\left\{\mean\left(\left.\sum_{j=1}^{\xi_i}\bone_{\{\ScS_i\le \eta_{ij}+1\}}\right|\xi_i,\ScS_i\right)\right\}\\
	&=&\sum_{i=1}^\infty\mean\left\{\xi_i\prob(\ScS_i\le \eta_{11}+1|\ScS_i)\right\}\\
	&=&\mean(\xi_1)\sum_{i=1}^\infty\sum_{l=0}^\infty \prob(\ScS_i\le l+1)q_l\\
	&=&r\mean(\zeta_1)\sum_{l=0}^\infty\mean\ScN_{l+1}q_l=r\mean(\zeta_1)\sum_{l=0}^\infty\frac{l+1}{\mean\tau_1}q_l\\
	&=&r\mean(\zeta_1)(\mean\eta_{11}+1)/\mean\tau_1.
	\end{eqnarray*}
	For $\V(\ScM)$, we write $\prob(\eta_{11}\wedge\eta_{12}=u)=q'_u$, then decomposing according to $i=l,\ j=m$; $i=l,\ j\ne m$; and $i\ne l$ for the second equation below, we obtain
	\begin{eqnarray}
	\mean(\ScM^2)&=&\mean\sum_{i=1}^\infty\sum_{l=1}^\infty\sum_{j=1}^{\xi_i}\sum_{m=1}^{\xi_l}\bone_{\{\ScS_i\le\eta_{ij}+1,\ScS_l\le\eta_{lm}+1\}}\nonumber\\
	&=&\mean\sum_{i=1}^\infty\sum_{j=1}^{\xi_i}\bone_{\{\ScS_i\le\eta_{ij}+1\}}+\mean\sum_{i=1}^\infty\sum_{1\le j\ne m\le\xi_i}\bone_{\{\ScS_i\le\eta_{ij}\wedge\eta_{im}+1\}}\nonumber\\
	&&+\mean\sum_{1\le i\ne l<\infty}\sum_{j=1}^{\xi_i}\sum_{m=1}^{\xi_l}\bone_{\{\ScS_i\le\eta_{ij}+1,\ScS_l\le\eta_{lm}+1\}}\nonumber\\
	&=&\mean\ScM+\sum_{i=1}^\infty\mean\left\{\mean\left(\left.\sum_{1\le j\ne m\le\xi_i}\bone_{\{\ScS_i\le\eta_{ij}\wedge\eta_{im}+1\}}\right|\ScS_i,\xi_i\right)\right\}\nonumber\\
	&&+\sum_{1\le i\ne l<\infty}\mean\left\{\mean\left(\left.\sum_{j=1}^{\xi_i}\sum_{m=1}^{\xi_l}\bone_{\{\ScS_i\le\eta_{ij}+1,\ScS_l\le\eta_{lm}+1\}}\right|\ScS_i,\ScS_l,\xi_i,\xi_l\right)\right\},\nonumber\end{eqnarray}
	which can be simplified to
	\begin{eqnarray}
	\mean(\ScM^2)&=&\mean\ScM+\sum_{i=1}^\infty\sum_{u=0}^\infty q'_u\mean\left\{\xi_i(\xi_i-1)\prob(\ScS_i\le u+1|\ScS_i)\right\}\nonumber\\
	&&+\sum_{1\le i\ne l<\infty}\sum_{u_1=0}^\infty\sum_{u_2=0}^\infty q_{u_1}q_{u_2}\mean\left\{\xi_i\xi_l\prob(\ScS_i\le u_1+1,\ScS_l\le u_2+1|\ScS_i,\ScS_l)\right\}\nonumber\\
	&=&\mean\ScM+\mean(\xi_1^2-\xi_1)\sum_{i=1}^\infty\sum_{u=0}^\infty q'_u\prob(\ScS_i\le u+1)\nonumber\\
	&&+(\mean\xi_1)^2\sum_{1\le i\ne l<\infty}\sum_{u_1=0}^\infty\sum_{u_2=0}^\infty q_{u_1}q_{u_2}\prob(\ScS_i\le u_1+1,\ScS_l\le u_2+1)\nonumber\\
	&=&\mean\ScM+\mean(\xi_1^2-\xi_1)\sum_{i=1}^\infty\sum_{u=0}^\infty q'_u\prob(\ScS_i\le u+1)\nonumber\\
	&&+(\mean\xi_1)^2\sum_{i=1}^\infty\sum_{l=1}^\infty\sum_{u_1=0}^\infty\sum_{u_2=0}^\infty q_{u_1}q_{u_2}\prob(\ScS_i\le u_1+1,\ScS_l\le u_2+1)\nonumber\\
	&&-(\mean\xi_1)^2\sum_{i=1}^\infty\sum_{u_1=0}^\infty\sum_{u_2=0}^\infty q_{u_1}q_{u_2}\prob(\ScS_i\le u_1\wedge u_2+1).\label{discretethm5}
	\end{eqnarray}
	Now, since $\mean\ScN_v=\sum_{i=1}^\infty\prob(\ScS_i\le v)$, $\mean\ScN_v=v/\mean\tau_1$, 
	$$\mean \ScN_{v_1}\ScN_{v_2}=\sum_{i=1}^\infty\sum_{l=1}^\infty\prob(\ScS_i\le v_1,\ScS_l\le v_2),$$
	we have
	$$\sum_{i=1}^\infty\sum_{u=0}^\infty q'_u\prob(\ScS_i\le u+1)
	=\sum_{u=0}^\infty q'_u\mean \ScN_{u+1}=\sum_{u=0}^\infty q'_u(u+1)/\mean\tau_1=[\mean(\eta_{11}\wedge\eta_{12})+1]/\mean\tau_1$$
	and
	\begin{eqnarray*}
	&&\sum_{i=1}^\infty\sum_{u_1=0}^\infty\sum_{u_2=0}^\infty q_{u_1}q_{u_2}\prob(\ScS_i\le u_1\wedge u_2+1)\\
	&&=\sum_{u_1=0}^\infty\sum_{u_2=0}^\infty q_{u_1}q_{u_2}\mean \ScN_{u_1\wedge u_2+1}\\
	&&=\sum_{u_1=0}^\infty\sum_{u_2=0}^\infty q_{u_1}q_{u_2}\frac{ u_1\wedge u_2+1}{\mean\tau_1}\\
	&&=[\mean(\eta_{11}\wedge\eta_{12})+1]/\mean\tau_1.
	\end{eqnarray*}
	 Thus, it follows from
	\Ref{discretethm5} that
	\begin{eqnarray}
	\mean(\ScM^2)
	&=&\mean\ScM+[\V(\xi_1)-\mean\xi_1][\mean(\eta_{11}\wedge\eta_{12})+1]/\mean\tau_1\nonumber\\
	&&+(\mean\xi_1)^2\sum_{u_1=0}^\infty\sum_{u_2=0}^\infty q_{u_1}q_{u_2}\mean \ScN_{u_1+1}\ScN_{u_2+1}.\label{discretethm6}
	\end{eqnarray}
	Define $\Delta \ScN_i=\ScN_i-\ScN_{i-1}$, $\Delta H(i)=H(i)-H(i-1)$. For $1\le m_1<m_2$, given $\Delta \ScN_{m_1}=1$, we know a renewal has happened at time $m_1$, so
	$$\mean(\Delta \ScN_{m_1}\Delta \ScN_{m_2})=\mean(\Delta \ScN_{m_2}|\Delta \ScN_{m_1}=1)\prob(\Delta \ScN_{m_1}=1)=\Delta H(m_2-m_1)/\mean\tau_1.$$
	This yields
	\begin{eqnarray}
	\mean \{\ScN_u^2\}&=&\mean\left\{\left(\sum_{m=1}^u\Delta\ScN_m\right)^2\right\}\nonumber\\
	&=&2\sum_{1\le m_1<m_2\le u}\mean\Delta\ScN_{m_1}\Delta\ScN_{m_2}+\sum_{m=1}^u\mean\{(\Delta\ScN_m)^2\}\nonumber\\
	&=&2\sum_{1\le m_1<m_2\le u}\Delta H(m_2-m_1)/\mean\tau_1+\sum_{m=1}^u\mean\{\Delta\ScN_m\}\nonumber\\
	&=&(\mean\tau_1)^{-1}\left(2\sum_{i=1}^{u-1}H(i)+u\right),\label{discretethm7}
	\end{eqnarray}
	and, for $u<v$,
	\begin{eqnarray}
	\mean(\ScN_u\ScN_v)&=&\mean(\ScN_u^2)+\mean[\ScN_u(\ScN_v-\ScN_u)]\nonumber\\
	&=&\mean(\ScN_u^2)+\sum_{m_1=1}^u\sum_{m_2=u+1}^v\mean(\Delta\ScN_{m_1}\Delta\ScN_{m_2})\nonumber\\
	&=&\mean(\ScN_u^2)+\sum_{m_1=1}^u\sum_{m_2=u+1}^v\Delta H(m_2-m_1)/\mean\tau_1\nonumber\\
	&=&(\mean\tau_1)^{-1}\left(\sum_{i=1}^{u-1}H(i)+u+\sum_{i=v-u}^{v-1}H(i)\right).\label{discretethm8}
	\end{eqnarray}
	Using \Ref{discretethm7} and \Ref{discretethm8}, we have
	\begin{eqnarray*}
	&&\sum_{v_1=0}^\infty\sum_{v_2=0}^\infty q_{v_1}q_{v_2}\mean \ScN_{v_1+1}\ScN_{v_2+1}\\
	&=&\sum_{v=0}^\infty q_v^2\mean (\ScN_{v+1}^2)+2\sum_{0\le v_1<v_2} q_{v_1}q_{v_2}\mean \ScN_{v_1+1}\ScN_{v_2+1}\\
	&=&(\mean\tau_1)^{-1}\sum_{v=0}^\infty\left(2\sum_{i=1}^vH(i)+v+1\right)q_v^2\\
	&&+
	2(\mean\tau_1)^{-1}\sum_{0\le v_1<v_2}\left(\sum_{i=1}^{v_1}H(i)+v_1+1+\sum_{i=v_2-v_1}^{v_2}H(i)\right)q_{v_1}q_{v_2}.
	\end{eqnarray*}
	This, together with \Ref{discretethm6}, \Ref{convrate05} and \Ref{convrate06}, gives \Ref{discretethm2}. For \Ref{discretethm3}, using the total probability formula, we have
	$$H(n)=\sum_{k=1}^n\mean(N_n|\tau_1=k)\prob(\tau_1=k)=\sum_{k=1}^n(H(n-k)+1)f_k=\sum_{k=1}^nH(n-k)f_k+F(n),$$
	as claimed. 
	\qed
	
	\section{Proof of Corollary 2}
	
	Let $\ScF$ be the $\sigma$-algebra generated by $\{\xi_i,\ i\ge 1\}$, then we have
	\begin{eqnarray}
	\phi_{\ScM}(u)&=&\mean\left\{\mean\left(\left.\prod_{i=1}^\infty\prod_{j=1}^{\xi_i}e^{\ime u\bonee_{\{\ScS_i\le \eta_{ij}+1\}}}\right|\ScF\right)\right\}\nonumber\\
	&=&\mean\prod_{i=1}^\infty\prod_{j=1}^{\xi_i}\left(e^{\ime u}(1-F_\eta(\ScS_i-2))+F_\eta(\ScS_i-2)\right)\nonumber\\
	&=&\mean\prod_{i=1}^\infty\left(e^{\ime u}+(1-e^{\ime u})F_\eta(\ScS_i-2)\right)^{\xi_i}\nonumber\\
	&=&\prod_{i=1}^\infty\psi_\xi\left(e^{\ime u}+(1-e^{\ime u})F_\eta(\ScS_i-2)\right),\label{cor2-3}
	\end{eqnarray}
	where $\psi_\xi(s)=\mean s^{\xi_1}$. However, $\xi_1=\sum_{j=1}^{\zeta_1}I_{1j}$ is a mixed binomial with
	\begin{equation}\psi_\xi(s)=\psi_\zeta(r(s-1)+1).\label{cor2-4}\end{equation}
	Hence, combining  \Ref{cor2-3} and \Ref{cor2-4} gives \Ref{cor2-1}. If $\zeta_i\sim\Pn(\lambda)$, we have
	$$\phi_{\ScM}(u)=e^{-r\lambda(1-e^{\imee u})\sum_{i=1}^\infty(1-F_\eta(\ScS_i-2))},$$
	which is the characteristic function of $\Pn\left(r\lambda\sum_{i=1}^\infty(1-F_\eta(\ScS_i-2))\right)$. \qed

	\section{Proof of Theorem 3}
	
	We have $\prob(\eta_{ij}\ge n+lD-i)=0$ for $n+lD-i>K$, hence
$$M_{n+lD} = \sum_{i=1}^{n+lD}\sum_{j=1}^{\xi_i}\bone_{\{\eta_{ij}\ge n+lD-i\}} =\sum_{i=n+lD-K}^{n+lD}\sum_{j=1}^{\xi_i}\bone_{\{\eta_{ij}\ge n+lD-i\}}\mbox{ almost surely}.$$
The characteristic function of $M_{n+lD}$ for $n+lD>K$ is
\begin{eqnarray*}
	\phi_{n+lD}(t)&=&\mean\left\{\mean\left(\left.\prod_{i=n+lD-K}^{n+lD}\prod_{j=1}^{\xi_i}e^{\ime t \bonee_{\{\eta_{ij} \geq n+lD-i \}}}\right|\xi_i,i=n+lD-K,\dots,n+lD \right)\right\} \\
	&=&\mean\prod_{i=n+lD-K}^{n+lD} \prod_{j=1}^{\xi_{i}}\left(e^{\ime t}(1-F_{\eta,i}(n+lD-i-1))+F_{\eta,i}(n+lD-i-1)\right)\\
	&=&\mean\prod_{i=n+lD-K}^{n+lD} \left(e^{\ime t}+(1-e^{\ime t})F_{\eta,i}(n+lD-i-1)\right)^{\xi_i}\\
	&=&\prod_{i=n+lD-K}^{n+lD}\psi_{\xi,i}\left(e^{\ime t}+(1-e^{\ime t})F_{\eta,i}(n+lD-i-1)\right)\\
	&=&\prod_{i=n+lD-K}^{n+lD}\psi_{\zeta,i}\left(r_i(e^{\ime t}-1)(1-F_{\eta,i}(n+lD-i-1))+1\right),
		\end{eqnarray*}
		where $\psi_{\xi,i}(s)=\mean s^{\xi_i}$ and the last equality is because the same reasoning as that for \Ref{cor2-4} gives 
		$$\psi_{\xi,i}(s)=\psi_{\zeta,i}(r_i(s-1)+1).$$

The mean of $M_{n+lD}$ is a straightforward calculation.  If $\zeta_i \sim \Pn(\lambda_i)$, the characteristic function $\phi_{n+lD}$ is reduced to
$$	\phi_{n+lD}(t) 
	= \exp \left( (e^{\ime t} - 1) \sum_{i=n+lD-K}^{n+lD} r_i\lambda_i(1-F_{\eta,i}(n+lD-i-1)) \right),
		$$
	the same as the characteristic function of $\Pn\left(\sum_{i=n+lD-K}^{n+lD}r_i\lambda_i(1-F_{\eta,i}(n+lD-i-1))\right)$.
\qed

\section{Plots}\label{app-plots}
  	
  	  		  	\begin{figure}[ht]
  	  \centering
  	    \includegraphics[width=.9\textwidth]{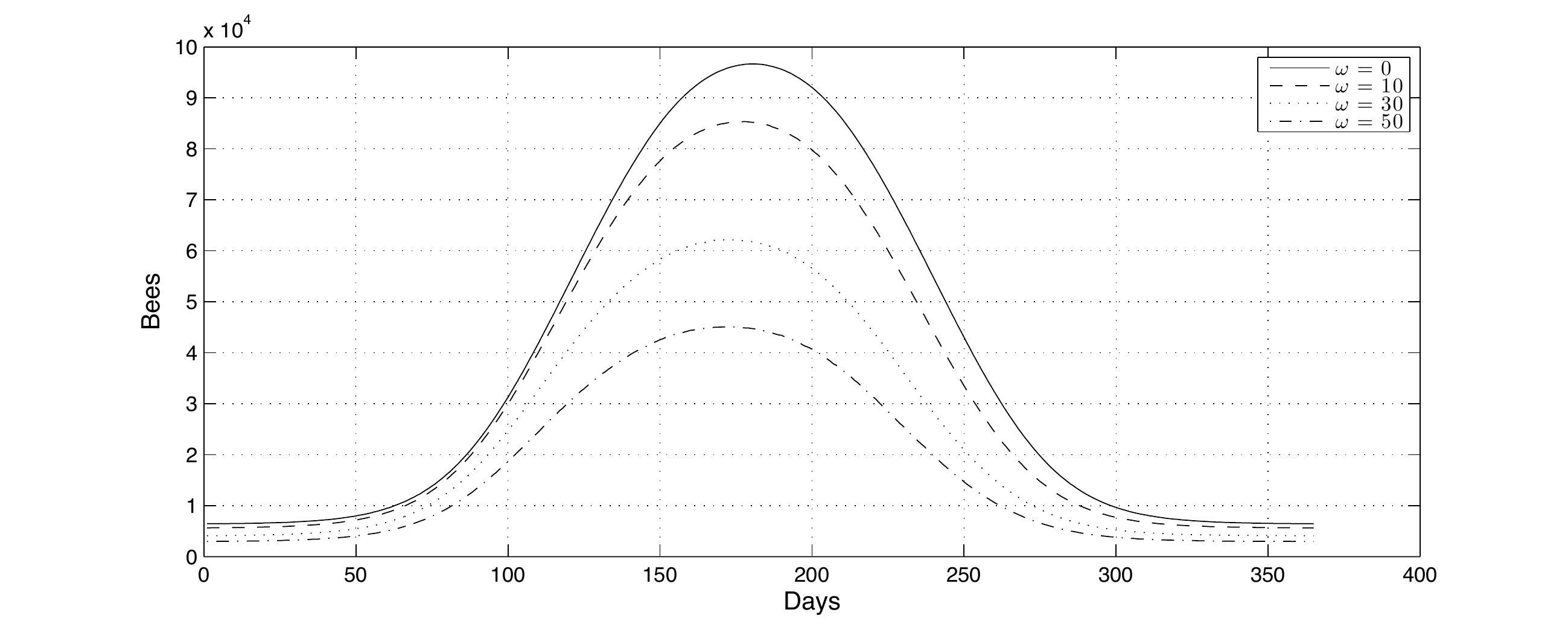}\\ a)\\
		 \includegraphics[width=.9\textwidth]{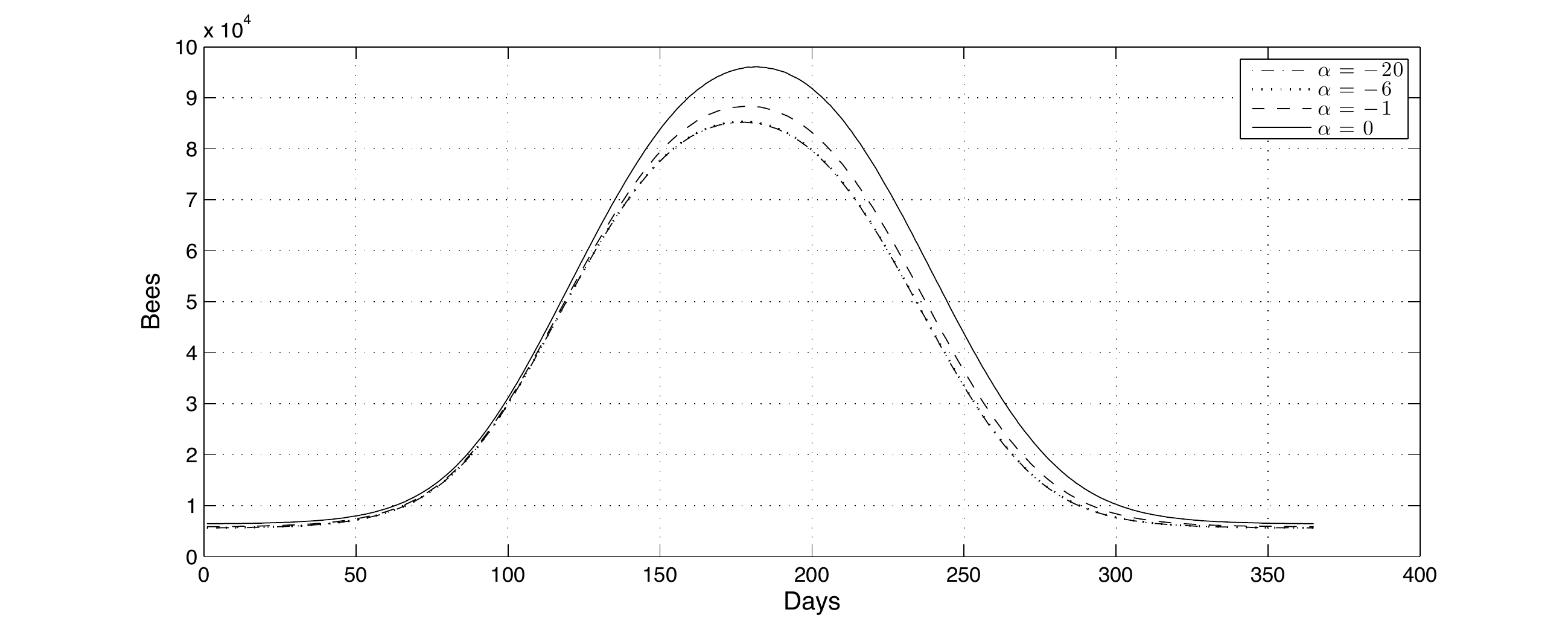}\\  b)\\
		   \includegraphics[width=.9\textwidth]{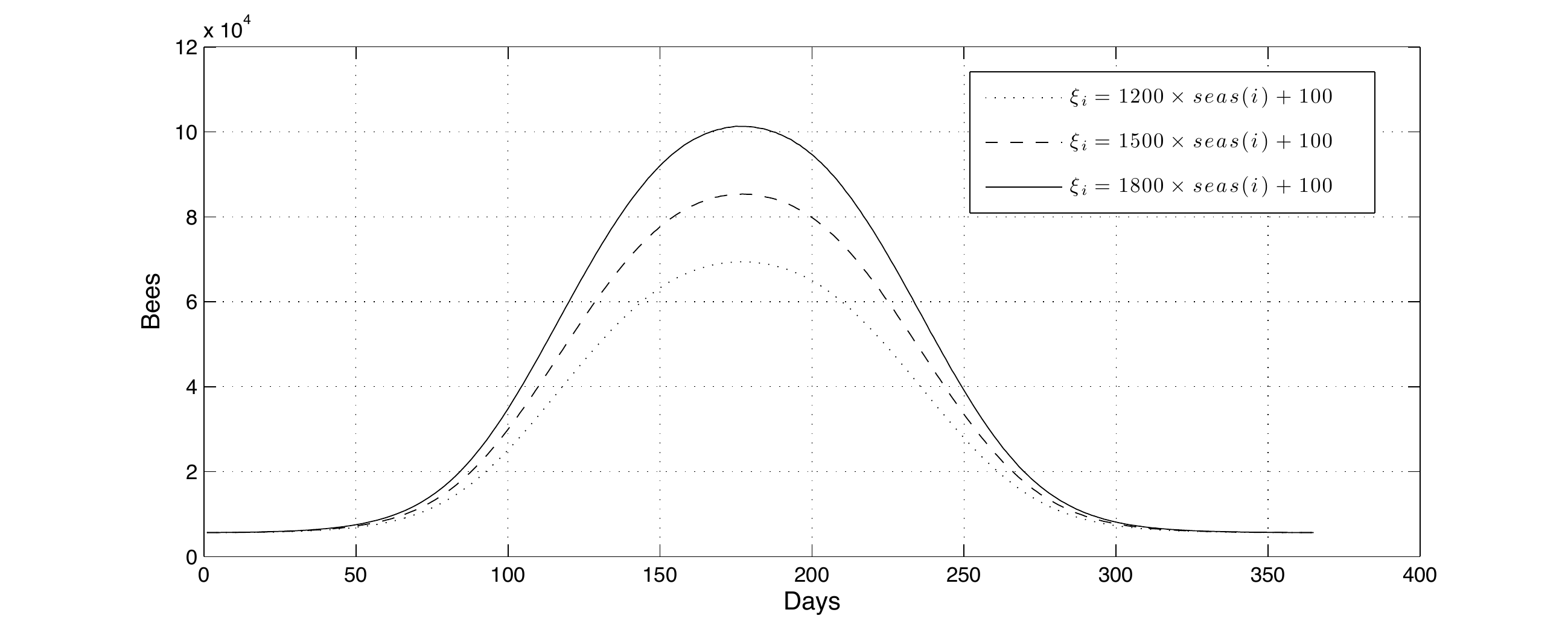}\\c)
  	  \caption{a) The mean of the daily colony size $M_n$ when  $\eta \sim SN(63,\omega^2,-6)$, for the scale parameter $\omega \in \{0,10,30,50\}$,
	  b) $\eta \sim SN(63,10^2,\alpha)$, for the slant parameter $\alpha \in \{-20,-6,-1,0\}$,  
	 and  a cyclic mean number of eggs: $\mean(\xi_i) {=r_i\mean(\zeta_i)}= 1500\times seas(i) + 100$.
	  c) $\eta \sim SN(63,10^2,-6)$, and a cyclic mean number of eggs given by $\mean(\xi_i) {=r_i\mean(\zeta_i)}= L\times seas(i) + 100, L\in\{1200, 1500, 1800\}$.}
  	  \label{fig:cyclic-mean-w}
  	\end{figure}

	\end{appendix}
		\end{document}